\documentclass{amsart}
\usepackage{amsthm}
\newtheorem{thm}{Theorem}[section]
\newtheorem{cor}[thm]{Corollary}
\newtheorem{lem}[thm]{Lemma}
\newtheorem{fact}[thm]{Fact}
\newtheorem{claim}[thm]{Claim}
\newtheorem{prop}[thm]{Proposition}

\newtheorem{defn}[thm]{Definition}

\numberwithin{equation}{section}

\begin{document}

\title[A second order smooth variational principle on Riemannian
manifolds]{A second order smooth variational principle on Riemannian manifolds}
\author{Daniel Azagra and Robb Fry}

\address{Departamento de An{\'a}lisis Matem{\'a}tico\\ Facultad de
Matem{\'a}ticas\\ Universidad Complutense\\ 28040 Madrid, Spain}

\address{Department of Mathematics and Statistics\\ School of Advanced
Technologies and Mathematics\\ Thompson Rivers University\\ 900
McGill Road Kamloops\\ British Columbia,  V2C 2N5 Canada}

\date{January 19, 2007}


\email{azagra@mat.ucm.es, rfry@tru.ca}



\begin{abstract}
We establish a second order smooth variational principle valid for
functions defined on (possibly infinite-dimensional) Riemannian
manifolds which are uniformly locally convex and have a strictly
positive injectivity radius and bounded sectional curvature.
\end{abstract}

\maketitle

\section{Introduction and main result}

It is well known that a continuous function defined on an
infinite-dimensional manifold (or on a Banach space) does not
generally attain a minimum in situations in which there would
typically exist minimizers if the function were defined on a
finite-dimensional manifold (for instance when the infimum of the
function in the interior of a ball is strictly smaller than the
infimum of the function on the boundary of the ball). In fact, as
shown in \cite{AC}, the smooth functions with no critical points
are dense in the space of continuous functions on every Hilbert
manifold (this result may be viewed as a strong approximate
version for infinite dimensional manifolds of the Morse-Sard
theorem). So, when we are given a smooth function on an
infinite-dimensional Riemannian manifold we should not expect to
be able to find any critical point, whatever the overall shape of
this function is, as there might be none.

This is quite inconvenient because many important problems of PDEs
and of optimization admit equivalent formulations as minimization
(or variational) problems, that is, one is given a continuous (or
differentiable, or convex, or lower semicontinuous, etc)
functional defined on a functional space, or on a (usually
infinite-dimensional) manifold, and one is asked to find a
minimizer of the functional, which will provide a solution of the
initial problem. Under several classes of rather restrictive
assumptions (for instance convexity of the functional and
reflexivity of the space) one may be able to show existence of
minimizers, but in many natural cases (such as the problem of
finding minimizing geodesics in infinite-dimensional Riemannian
manifolds, to name one) exact minimizers do not generally exist,
and one is thus forced to look for {\em approximate minimizers},
which will provide {\em approximate solutions} of the problem.

This is why perturbed minimization principles, or variational
principles, are important. A variational principle typically
asserts that, for a given lower semicontinuous function
$f:\mathcal{X}\to\mathbb{R}$, defined on a given space
$\mathcal{X}$, and such that $f$ is bounded below, there exists a
function $\varphi:\mathcal{X}\to\mathbb{R}$ belonging to a given
class $\mathcal{C}$ such that $f-\varphi$ attains a global minimum
and $\varphi$ can be prescribed to be small in some sense.
Historically, Ekeland's variational principle \cite{Ekeland} was
the first perturbed minimization principle discovered, and one of
the most powerful and widely applicable, because it holds for any
complete metric space $\mathcal{X}$; in this case the perturbing
class $\mathcal{C}$ is the set of small scalar multiples of
distance functions to points, that is the perturbing function
$\varphi$ has the shape of a flat cone. One of the many striking
applications of Ekeland's principle was to show the existence of
almost minimizing geodesics on complete infinite-dimensional
Riemannian manifolds, that is an approximate version of the
Hopf-Rinow theorem (as we recalled above, an exact analogue is
false in infinite dimensions).

Despite the generality and power of Ekeland's result, there are
some limitations to the applicability of this minimization
principle due to the fact that the perturbing functions $\varphi$
are not differentiable. Indeed, there are many situations (for
instance when one wants to build a theory of subdifferentiability,
in nonsmooth analysis, and in the study of viscosity solutions to
Hamilton-Jacobi equations) in which one needs the perturbing
function $\varphi$ to be differentiable. In order to remedy this
deficiency Borwein and Preiss \cite{BorweinPreiss} established a
{\em smooth} variational principle in which the space
$\mathcal{X}$ is a Banach space with a $C^1$ smooth norm, and the
perturbing functions $\varphi$ are smooth functions which can be
taken to be arbitrarily small and with an  arbitrarily small
Lipschitz constant. Later on, Deville, Godefroy and Zizler
\cite{DGZ, DGZsvp}, by using a new method of proof based on the
use of Baire's category theorem, were able to extend this smooth
variational principle to the class of all Banach spaces with
smooth bump functions, as well as for higher orders of smoothness.

Smooth variational principles on Riemannian manifolds were not
studied until very recently. In \cite{AFL2} a Riemannian version
of the Deville-Godefroy-Zizler smooth variational principle was
established within the class of complete Riemannian manifolds
which are {\em uniformly bumpable}. This was applied to developing
a theory of nonsmooth analysis and to the study of viscosity
solutions to Hamilton-Jacobi equations on Riemannian manifolds.
Later, in \cite{AFLR}, the authors showed that the assumption on
uniform bumpability can be dispensed with, by establishing that
every Lipschitz function $f$ defined on a Riemannian manifold can
be uniformly approximated by a sequence of $C^\infty$ smooth
Lipschitz functions $(f_{n})$ in such a way that the Lipschitz
constants of $f_{n}$ converge to the Lipschitz constant of $f$.

\medskip

In this paper we will prove a second order smooth variational
principle valid for functions defined on every complete Riemannian
manifold of bounded sectional curvature, uniformly locally convex,
and with a positive injectivity radius. In fact the result is
valid (see Theorem \ref{soub manifolds satisfy the sosvp} below)
for every complete Riemannian manifold which is {\em second order
uniformly bumpable} according to Definition \ref{definition of
uniform bumpability} given below. The question whether or not
every complete Riemannian manifold is second order uniformly
bumpable is open, and seems to be much more subtle than its first
order analogue.

Let us state our main result.

\begin{thm}\label{main theorem}
Let $M$ be a (possibly infinite-dimensional) complete Riemannian
manifold. Assume that the sectional curvature of $M$ is bounded,
and that $M$ is uniformly locally convex and has a positive
injectivity radius. Then, for every lower semicontinuous function
$f:M\to (-\infty, \infty]$ which is bounded below, with
$f\not\equiv+\infty$, and for every $\varepsilon>0$, there exists
a $C^2$ smooth function $\varphi: M\to\mathbb{R}$ such that
\begin{enumerate}
\item $f-\varphi$ attains its strong minimum on $M$;
\item $\|\varphi\|_{\infty}<\varepsilon$;
\item $\|d\varphi\|_{\infty}<\varepsilon$;
\item $\|d^{2}\varphi\|_{\infty}<\varepsilon$.
\end{enumerate}
\end{thm}
\noindent Recall that a function $g$ is said to attain a strong
minimum at a point $x_{0}$ provided $g(x_{0})=\inf_{x\in M}g(x)$
and $\lim_{n\to\infty}d(x_{n},x_{0})=0$ whenever $(x_{n})$ is a
minimizing sequence (that is, if \,
$\lim_{n\to\infty}g(x_{n})=g(x_{0})$).

In the statement of Theorem \ref{main theorem} we used the
following notation:
\begin{eqnarray*}
& & \|\varphi\|_{\infty}=\sup_{x\in M}|\varphi(x)| \\
& & \|d\varphi\|_{\infty}=\sup_{x\in M}\|d\varphi(x)\|_{x} \\
& & \|d^{2}\varphi\|_{\infty}=\sup_{x\in
M}\|d^{2}\varphi(x)\|_{x},
\end{eqnarray*}
where
    $$
    \|d\varphi(x)\|=\sup_{v\in TM_{x}, \|v\|_{x}=1}d\varphi(x)(v)=\sup_{v\in
    TM_{x}, \|v\|_{x}=1}\langle\nabla\varphi(x), v\rangle_{x},
    $$
and
    $$
    \|d^{2}\varphi(x)\|_{x}=\sup_{v\in TM_{x},
    \|v\|_{x}=1}|d^{2}\varphi(x)(v,v)| .
    $$
Recall also that the Hessian $D^{2}\varphi$ of a $C^2$ smooth
function $\varphi$ on $M$ is defined by
    $$
    D^{2}\varphi(X,Y)=\langle\nabla_{X}\nabla\varphi, Y\rangle,
    $$
where $\nabla\varphi$ is the gradient of $\varphi$; $X$, $Y$ are
vector fields on $M$, and $\nabla_{X}Z$ is the covariant
derivative of a vector field $Z$ with respect to a vector field
$X$. The Hessian is a symmetric tensor field of type $(0,2)$ and,
for a point $p\in M$, the value $D^{2}\varphi(X,Y)(p)$ only
depends of $f$ and the vectors $X(p), Y(p)\in TM_{p}$. So we can
define the second derivative of $\varphi$ at $p$ as the symmetric
bilinear form $d^{2}\varphi(p):TM_{p}\times TM_{p}\to\mathbb{R}$
    $$
    (v,w)\mapsto d^{2}\varphi(p)(v,w):=D^{2}\varphi(X,Y)(p),
    $$
where $X, Y$ are any vector fields such that $X(p)=v, Y(p)=w$. A
very useful way to compute $d^{2}\varphi(p)(v,v)$ is to take a
geodesic $\gamma$ with $\gamma'(0)=v$ and calculate
    $$
    \frac{d^{2}}{dt^{2}} \varphi(\gamma(t))|_{t=0},
    $$
which equals $d^{2}\varphi(p)(v,v)$. We will often write
$d^{2}\varphi(p)(v)^{2}$ instead of $d^{2}\varphi(p)(v,v)$.

\medskip

Let us fix some notation and terminology that we will be using
throughout the paper. $M=(M, g)$ will always be a (possibly
infinite-dimensional) complete Riemannian manifold, and
$g(x)=\langle \cdot, \cdot\rangle_{x}$ its Riemannian metric. We
refer to \cite{Lang} for a basic introduction to
infinite-dimensional Riemannian manifolds. However, instead of
Lang's notation we will employ a more classic notation, such as
that of Do Carmo's book \cite{dC}. The letters $X, Y, Z, V, W$
will stand for smooth vector fields on $M$, and $\nabla_{Y}X$ will
always denote the covariant derivative of $X$ with respect to $Y$.
The Riemannian curvature of $M$ will be denoted by $R$. Recall
that the value of
$R(X,Y)Z:=\nabla_{X}\nabla_{Y}Z-\nabla_{Y}\nabla_{X}Z-\nabla_{[X,Y]}Z$
at a point $p\in M$ only depends on the values of $X, Y, Z$ at
$p$, and that the sectional curvature of $M$ at a point $p$ with
respect a plane spanned by two vectors $v, w\in TM_{p}$ is defined
by
    $$
    K(p; v,w)=\frac{\langle R(v,w)v, v \rangle}{|v\wedge w|^{2}}
    $$
where $|v\wedge w|$ is the area of the parallelogram defined by
$u, v$ in $TM_{p}$.

Geodesics in $M$ will be denoted by $\gamma$, $\sigma$, and their
velocity fields by $\gamma', \sigma'$. If $X$ is a vector field
along $\gamma$ we will often denote $X'(t)=
\frac{D}{dt}X(t)=\nabla_{\gamma'(t)}X(t)$. Recall that $X$ is said
to be parallel along $\gamma$ if $X'(t)=0$ for all $t$. The
Riemannian distance in $M$ will always be denoted by $d(x,y)$
(defined as the infimum of the lengths of all curves joining $x$
to $y$ in $M$).

We will often identify the tangent space of $M$ at a point $x$,
denoted by $TM_{x}$, with the cotangent space at $x$, denoted by
$TM_{x}^{*}$. The space of bilinear forms on $TM_{x}$
(respectively symmetric bilinear forms) will be denoted by
$\mathcal{L}^{2}(TM_{x})$ or $\mathcal{L}^{2}(TM_{x}, \mathbb{R})$
(resp. $\mathcal{L}^{2}_{s}(TM_{x})$ or
$\mathcal{L}^{2}_{s}(TM_{x}, \mathbb{R})$). Also, we will denote
by $T_{2, s}(M)$ the tensor bundle of symmetric bilinear forms,
that is
$$
T_{2, s}(M)=\bigcup_{x\in M}\mathcal{L}^{2}_{s}(TM_x, \mathbb{R}),
$$
and $T_{2, s}(M)_x = \mathcal{L}^{2}_{s}(TM_x, \mathbb{R})$.

We will make extensive use of the exponential mapping $\exp_{x}$
throughout the paper. Recall that for every $x\in M$ there exists
a mapping $\exp_{x}$, defined on a neighborhood of $0$ in the
tangent space $TM_x$, and taking values in $M$, which is a local
diffeomorphism and maps straight line segments passing through $0$
onto geodesic segments in $M$ passing through $x$.

By $i_{M}(x)$ we will denote the {\em injectivity radius of $M$ at
$x$}, that is the supremum of the radius $r$ of all balls
$B(0_{x}, r)$ in $TM_{x}$ for which $\exp_{x}$ is a diffeomorphism
from $B(0_{x}, r)$ onto $B(x,r)$. The number $i_{M}(x)$ is
strictly positive for every $x$. Similarly, $i(M)$ will denote the
{\em global injectivity radius of $M$}, specifically
$i(M)=\inf\{i_{M}(x)
: x\in M\}$.

We will also need to recall some results about convexity in
Riemannian manifolds. We say that a subset $U$ of a Riemannian
manifold is {\em convex} if given $x, y\in U$ there exists a
unique geodesic in $U$ joining $x$ to $y$, and such that the
length of the geodesic is $d(x,y)$. Every Riemannian manifold is
{\em locally convex}, in the sense that for every $x\in M$, there
exists $c>0$ such that for all $r$ with $0<r<c$, the open ball
$B(x,r)=\exp_{x} B(0_{x}, r)$ is convex (this is Whitehead's
Theorem). The {\em convexity radius} of a point $x\in M$ in a
Riemannian manifold $M$ is defined as the supremum in
$\overline{\mathbb{R}^{+}}$ of the numbers $r>0$ such that the
ball $B(x,r)$ is convex. We denote this supremum by $c_{M}(x)$ (by
the result we have just mentioned, $c_{M}(x)$ is strictly positive
for every $x$).

Whitehead's theorem gives rise to the notion of {\em uniformly
locally convex} manifold: we say that a Riemannian manifold $M$ is
{\em uniformly locally convex} provided that there exists $c>0$
such that for every $x\in M$ and every $r$ with $0<r<c$ the ball
$B(x,r)=\exp_{x} B(0_{x}, r)$ is convex. This amounts to saying
that the {\em global convexity radius} of $M$ (defined as
$c(M):=\inf\{ c_{M}(x) : x\in M\}$) is strictly positive.

It is also worth noting that the functions $x\mapsto i_{M}(x)$ and
$x\mapsto c_{M}(x)$ are continuous, see \cite{Klingenberg}.

As noted above, the Hopf-Rinow theorem fails in infinite
dimensions, and hence generally geodesics minimizing the distance
between two given points do not exist in infinite-dimensional
Riemannian manifolds. But, if $M$ is locally uniformly convex and
$i(M)>0$, they do always exist locally. In fact they exist in a
{\em uniformly local} way: for any $r>0$ with $r<\min\{i(M),
c(M)\}$ and for all points $x, y$ with $d(x,y)\leq r$, there is a
unique (up to reparameterizations) minimizing geodesic $\gamma$
connecting $x$ and $y$ ($\gamma$ is defined by
$\gamma(t)=\exp_{x}(tw_{y})$, where $w_{y}=\exp_{x}^{-1}(y)$).
Moreover, $d(x,y)$ is given by
$d(x,y)^{2}=\|\exp_{x}^{-1}(y)\|_{x}$, see \cite[Proposition
3.9]{AFL2}. This will be extensively used in Section $2$ of the
paper.

Let us also recall that for a minimizing geodesic $\gamma:[0,
\ell]\to M$ connecting $x$ to $y$ in $M$, and for a vector $v\in
TM_{x}$ there is a unique parallel vector field $P$ along $\gamma$
such that $P(0)=v$. The mapping $TM_{x}\ni v\mapsto P(\ell)\in
TM_{y}$ is called the {\em parallel translation} of $v$ along
$\gamma$ and is a linear isometry from $TM_{x}$ onto $TM_{y}$
which we will denote by $L_{xy}$.

For a basic theory of Jacobi fields on infinite-dimensional
Riemannian manifolds (and for any other unexplained terms of
Riemannian geometry used in Section $2$), we refer the reader to
\cite{Lang}.

\medskip

Let us finish this introduction by providing a typical application
of the smooth variational principle, namely to show the existence
and density of the set of points of subdifferentiability (of
second order, in our case). For a lower semicontinuous function
$f:M\to (-\infty, +\infty]$ we define the second order subjet of
$f$ at a point $x\in M$ as the set
    $$
    J^{2,-}f(x)=\{(d\varphi(x), d^{2}\varphi(x)) \, : \,
    \varphi\in C^{2}(M, \mathbb{R}), \,  f-\varphi \textrm{
    attains a local minimum at } x\}.
    $$
\begin{cor}
Let $M$ be as in the statement of Theorem \ref{main theorem}, and
let $f:M\to (-\infty, +\infty]$ be a lower semicontinuous
function. Then the set $\{z\in M
: J^{2,-}f(z)\neq\emptyset\}$ is dense in the set $\{x\in M:
f(x)<\infty\}$.
\end{cor}
\begin{proof}
The proof is an easy consequence of Theorem \ref{main theorem}.
Indeed, the same proof as in\cite[Proposition 4.17]{AFL2} holds if
we replace $D^{-}f(x)$ with $J^{2,-}f(x)$, and $C^1$ with $C^2$.
\end{proof}

\medskip

\section{Second order uniformly bumpable manifolds}

In this section we will show that all complete Riemannian
manifolds which are uniformly locally convex and have bounded
sectional curvature and a strictly positive injectivity radius are
{\em second order uniformly bumpable}, meaning the following.

\begin{defn}\label{definition of uniform bumpability}
{\em We will say that a Riemannian manifold
$M$ is {\em second order uniformly bumpable} provided there exist
$R>1, r>0$ such that for every $z\in M$ and for every $\delta\in
(0, r)$ there exists a $C^2$ smooth function $b:M\to [0,1]$ such
that
\begin{enumerate}
\item[{(i)}] $b(z)=1$
\item[{(ii)}] $b(x)=0$ whenever $d(x,z)\geq \delta$
\item[{(iii)}] $\|db\|_{\infty}\leq R/\delta$
\item[{(iv)}] $\|d^{2}b\|_{\infty}\leq R/\delta^{2}$.
\end{enumerate}
}
\end{defn}

The key to the proof of the main result of this section is the
following estimation for the norm of the second derivative of the
square of the distance function, for which we have found no
reference (of course there are well known estimations in the
finite-dimensional case, see for instance \cite[Lemma 2.9 and
Exercise 4, page 153]{Sakai}, but they all depend on results
established exclusively in finite dimensions, and of which we have
found no infinite-dimensional versions in the literature).

\begin{prop}\label{estimation for the second derivative of the squared distance}
If $M$ is a Riemannian manifold whose sectional curvature $K$ is
bounded, say $|K|\leq K_{0}$, and $c(M)>0$, $i(M)>0$, then, for
every $r$ with $0<r<\min\{i(M), c(M), \pi/2\sqrt{K_{0}}\, \}$, and
for every $z\in M$, the function $\varphi(x):= d(x,z)^{2}$ is
$C^\infty$ smooth on $B(z, r)$ and its second derivative satisfies
    $$
\|d^{2}\varphi(x)\|_{x}\leq (2+\frac{2}{3}K_{0}\, d(x,z)^{2})
    $$
for all $x\in B(z,r)$.
\end{prop}
In the proof of this Proposition we will have to use some well
known results about the second variation of the arc length and the
energy functionals. Let us briefly review the facts that we will
be using.

Fix a number $r$ such that $0<r<\min\{i(M), c(M),
\pi/2\sqrt{K_{0}}\,\}$, and take two points $x, x_{0}\in M$ with
$d(x, x_{0})<r$. Let $\gamma$ be the unique minimizing geodesic,
parameterized by arc-length, connecting $x_{0}$ to $x$. Denote
$\ell=d(x, x_{0})$, the length of $\gamma$. Consider
$\alpha(t,s)$, a smooth variation of $\gamma$, that is a smooth
mapping $\alpha:[0,\ell]\times [-\varepsilon, \varepsilon]\to M$
such that $\alpha(t,0)=\gamma(t)$ for all $t\in [0,\ell]$.
Consider the length and the energy functionals, defined by
    $$
    L(s)=L(\alpha_{s})=\int_{0}^{\ell}\|\alpha_{s}'(t)\|dt
    $$
and
    $$
    E(s)=E(\alpha_{s})=\int_{0}^{\ell}\|\alpha_{s}'(t)\|^{2}dt,
    $$
where $\alpha_{s}$ is the variation curve defined by
$\alpha_{s}(t)=\alpha(t,s)$ for every $t\in [0,\ell]$. According
to the Cauchy-Schwarz inequality (applied to the functions
$f\equiv 1$ and $g(t)=\|\alpha_{s}'(t)\|$ on the interval $[0,
\ell ]$) we have that
    $$
    L(s)^{2}\leq\ell E(s), \eqno(0)
    $$
with equality if and only if $\|\alpha_{s}'(t)\|$ is constant. In
particular we have $$L(0)^{2}=\ell E(0)$$ because
$\alpha_{0}=\gamma$ is a geodesic.

Therefore, if we furthermore assume that $\alpha_{s}$ is a
geodesic for each $s$ (that is $\alpha$ is a variation of $\gamma$
through geodesics) we have that
    $$
    L(s)^{2}=\ell E(s) \eqno(1)
    $$
for every $s\in [-\varepsilon, \varepsilon]$.

Now define $\varphi(x)=d(x,x_{0})^{2}$. The function $\varphi$ is
$C^\infty$ smooth on $B(x_{0}, r)$ because $r<\min\{i(M), c(M)\}$,
so on this ball $\exp_{x_{0}}^{-1}$ is a $C^\infty$
diffeomorphism, and
    $$
    d(x,x_{0})^{2}=\|\exp_{x_{0}}^{-1}(x)\|_{x_{0}}^{2}.
    $$
Let us take a vector $v\in TM_{x}$. We want to calculate
    $$
    d^{2}\varphi(x)(v)^{2},
    $$
which is given by
    $$
    \frac{d^{2}}{ds^{2}}\varphi(\sigma_{x}(s))|_{s=0},
    $$
where $\sigma_{x}(s)=\exp_{x}(sv)$. To this end let us denote by
$\alpha_{s}:[0,\ell]\to M$ the unique minimizing geodesic joining
the point $x_{0}$ to the point $\sigma_{x}(s)$ (notice that now,
for $s\neq 0$, $\alpha_{s}$ is not necessarily parameterized by
arc-length), and let us define $\alpha:[0,\ell]\times
[-\varepsilon, \varepsilon]\to M$ by $\alpha(t,s)=\alpha_{s}(t)$.
Then $\alpha$ is a smooth variation through geodesics of
$\gamma(t)=\alpha(t,0)$ and we have
    $$
    \varphi(\sigma_{x}(s))=L(s)^{2}=\ell
    E(s),
    $$
and therefore
    $$
    d^{2}\varphi(x)(v)^{2}=\ell E''(0). \eqno(2)
    $$
If we denote $X(t)=\partial\alpha(t,0)/\partial s$, the
variational field of $\alpha$, then the formula for the second
variation of energy (see \cite{Lang}, Chapter XI) tells us that
    $$
    \frac{1}{2}E''(0)=\int_{0}^{\ell}\left( \langle X', X'\rangle -
    \langle R(\gamma',X)\gamma', X\rangle \right)dt \, + \,
    \langle\frac{D}{ds}\frac{\partial\alpha}{\partial s}(t,0),
    \gamma'(t)\rangle |_{t=0}^{t=\ell}, \eqno(3)
    $$
or equivalently
    $$
    \frac{1}{2}E''(0)=-\int_{0}^{\ell}\langle X,
    X''+R(\gamma',X)\gamma'\rangle dt \, + \, \langle X(t),
    X'(t)\rangle |_{t=0}^{t=\ell} \, + \,
    \langle\frac{D}{ds}\frac{\partial\alpha}{\partial s}(t,0),
    \gamma'(t)\rangle |_{t=0}^{t=\ell}, \eqno(3')
    $$
where we denote $X'=\nabla_{\gamma'(t)}X$, and
$X''=\nabla_{\gamma'(t)} X'$.

Note that, since the variation field of a variation through
geodesics is always a Jacobi field, and since the points $x$ and
$x_{0}$ are not conjugate (recall that $r<i(M)$), $X$ is the
unique Jacobi field along $\gamma$ satisfying that $X(0)=0$,
$X(\ell)=v$, that is $X$ is the unique vector field along $\gamma$
satisfying
    $$
    X''(t)+R(\gamma'(t), X(t))\gamma'(t)=0, \,\,\, \textrm{ and } \,\,\, X(0)=0, \,\,\
    X(\ell)=v,
    $$
where $R$ is the curvature of $M$.

On the other hand, since the curve $s\to\alpha(\ell,
s)=\sigma_{x}(s)$ and $s\to\alpha(0, s)\equiv x_{0}$ are
geodesics, we have that
    $$
    \langle\frac{D}{ds}\frac{\partial\alpha}{\partial s}(t,0),
    \gamma'(t)\rangle |_{t=0}^{t=\ell}=0.
    $$
These observations allow us to simplify the formulas $(3)$ and
$(3')$ by dropping the terms that vanish, thus obtaining that
    $$
    \frac{1}{2}E''(0)=\int_{0}^{\ell}\left( \langle X', X'\rangle -
    \langle R(\gamma',X)\gamma', X\rangle \right)dt, \eqno(4)
    $$
and also
    $$
    \frac{1}{2}E''(0)=\langle X(\ell),
    X'(\ell)\rangle. \eqno(4')
    $$
Recall that the right-hand side of $(4)$ is called the index form
and is denoted by $I(X,X)$.

By combining $(2)$ and $(4)$ we get
    $$
    d^{2}\varphi(x_{0})(v)^{2}=2\ell \int_{0}^{\ell}\left( \langle X', X'\rangle -
    \langle R(\gamma',X)\gamma', X\rangle \right)dt, \eqno(5)
    $$
or equivalently
    $$
    d^{2}\varphi(x_{0})(v)^{2}=2\ell\langle X'(\ell), X(\ell)\rangle. \eqno(5')
    $$
We are going to use these formulas to deduce the estimation in the
statement of Proposition \ref{estimation for the second derivative
of the squared distance}. But we will need to combine them with a
couple of facts about Jacobi fields, and with the Rauch Comparison
Theorem. First, we must use the following.
\begin{lem}\label{the index form is positive definite for minimizing geodesics}
Let $\gamma:[0,\ell]\to M$ be a geodesic whose length, $\ell$, is
the distance between its end points, and let $Y$ be a field along
$\gamma$ with $Y(0)=0$ and $Y(\ell)=0$. Then
    $$
    I(Y,Y)\geq 0.
    $$
\end{lem}
\begin{proof}
This result is stated and proved in \cite[Theorem 1.7 of Ch.
XI]{Lang} under the additional assumption that $Y$ is orthogonal
to $\gamma$, but it is true for any $Y$ (and in fact almost the
same proof holds, with an additional remark). We will write the
whole proof for the reader's convenience. Define
    $$
    \beta(t,s)=\exp_{\gamma(t)}(s Y(t))
    $$
for $0\leq t\leq\ell$ and $0\leq s\leq\varepsilon$, and for a
sufficiently small $\varepsilon>0$. For each $s$,
$\beta_{s}(t):=\beta(t,s)$ is a curve, not necessarily a geodesic,
joining the end points of $\gamma$, that is
    $$
    \beta_{s}(0)=\gamma(0), \textrm{ and } \beta_{s}(\ell)=\gamma(\ell)
    $$
(because of the assumption that $Y(0)=0, Y(\ell)=0$). Moreover,
$\beta(t,0)=\gamma(t)$, so $\beta$ is a variation of $\gamma$
which leaves the end points fixed, and the variation field of
$\beta$ is
    $$
    \frac{\partial\beta(t,0)}{\partial s}=Y(t).
    $$
Using the formula for the second variation of energy \cite[Chapter
XI]{Lang} (and taking into account that the curves
$s\mapsto\beta(0,s)$ and $s\mapsto\beta(\ell,s)$ are geodesics, in
fact points), we have that
    $$
    \frac{1}{2}E''(0)=\int_{0}^{\ell}\left( \langle Y', Y'\rangle -
    \langle R(\gamma',Y)\gamma', Y\rangle \right)dt:=I(Y,Y).
    \eqno(6)
    $$
On the other hand, according to equation $(0)$, and bearing in
mind that $\gamma$ is a minimizing geodesic, we have that
    $$
    \ell E(0)=L(0)^{2}\leq L(s)^{2}\leq\ell E(s),
    $$
hence
    $$
    E(0)\leq E(s).
    $$
That is, a geodesic which minimizes length also minimizes energy.
Therefore we have
    $$
    E''(0)\geq 0,
    $$
which, combined with equation $(6)$ yields $I(Y,Y)\geq 0$.
\end{proof}
This Lemma allows us to improve Corollary 1.8 of \cite[Chapter
XI]{Lang} as follows.
\begin{prop}\label{Spivaks lemma}
Let $\gamma:[0,\ell]\to M$ be a geodesic whose length, $\ell$, is
the distance between its end points. Let $X$ be a Jacobi field
along $\gamma$, and $Z$ any smooth vector field along $\gamma$
such that $X(0)=Z(0)$ and $X(\ell)=Z(\ell)$. Then
    $$
    I(X, X)\leq I(Z,Z).
    $$

\noindent {\em In summary, among all vector fields along $\gamma$
with the same boundary conditions, the unique Jacobi field along
$\gamma$ determined by those conditions minimizes the index form.}
\end{prop}
\begin{proof}
The same proof as in \cite[Chapter XI, Corollary 1.8]{Lang} holds
(the additional assumption in \cite{Lang} that $X-Z$ is orthogonal
to $\gamma$ is only made to be able to apply Theorem 1.7 in
Chapter XI of \cite{Lang}, which we have just improved in the
above Lemma by dispensing with that extra assumption). See also
Corollary 10 in Chapter 8 of \cite{Spivak}, where the
finite-dimensional version of this Proposition is established with
the same proof (which only uses the finite-dimensional version of
the above Lemma and the second variation formula).
\end{proof}

As noted above, we will also need to use an infinite-dimensional
version of the Rauch Comparison Theorem, which we next state for
the reader's convenience.
\begin{thm}[Rauch Comparison Theorem]\label{Rauch}
Let $M$ and $\widetilde{M}$ be Riemannian manifolds of the same
dimension, which may be infinite. Let $\gamma$ (resp.
$\widetilde{\gamma}$) be geodesics in $M$ (resp. $\widetilde{M}$),
parameterized by arc-length, and defined on the same interval $[0,
\ell]$. Let $J$ (resp $\widetilde{J}$) be Jacobi fields along
these geodesics, orthogonal to $\gamma$ (resp.
$\widetilde{\gamma}$). Assume that
\begin{enumerate}
\item
$J(0)=\widetilde{J}(0)=0$, and $J(t), \widetilde{J}(t)\neq 0$ for
$0<t\leq\ell$.
\item $\|J'(0)\|=\|\widetilde{J}'(0)\|$.
\item The length of $\gamma$ is the distance between its end
points.
\item The sectional curvature of $M$ is less than or equal to
the sectional curvature of $\widetilde{M}$.
\end{enumerate}
Then, for all $t\in (0, \ell]$, we have
    $$
    \|\widetilde{J}(t)\|\leq \|J(t)\|,
    $$
and
    $$
    \frac{\langle \widetilde{J}'(t), \widetilde{J}(t)\rangle}{\langle \widetilde{J}(t), \widetilde{J}(t)\rangle}
    \leq\frac{\langle J'(t), J(t)\rangle}{\langle J(t), J(t)\rangle}.
    $$
\end{thm}
\noindent A proof of this result can be found in \cite[p.
319]{Lang}. The last inequality does not appear in Lang's
statement of \cite[Chapter XI, Theorem 5.1]{Lang}, but is
established and used in his proof.

Now we are ready to provide a

\begin{center}
{\bf Proof of Proposition \ref{estimation for the second
derivative of the squared distance}.}
\end{center}

\noindent According to equations $(5)$ and $(5')$
    $$
    d^{2}\varphi(x_{0})(v)^{2}=2\ell \int_{0}^{\ell}\left( \langle X', X'\rangle -
    \langle R(\gamma',X)\gamma', X\rangle \right)dt \, = 2\ell\, \langle X'(\ell), X(\ell)\rangle,
    $$
where $X$ is the unique Jacobi field along $\gamma$ such that
$X(0)=0$ and $X(\ell)=v$ (and $\gamma:[0, \ell]\to M$ is the
unique minimizing geodesic, parameterized by arc length, such that
$\gamma(0)=x_{0}, \gamma(\ell)=x$). Recall that we are assuming
$d(x,x_{0})\leq r<\min\{i(M), c(M), \pi/2\sqrt{K_{0}}\, \}$).

With the help of Proposition \ref{Spivaks lemma} and the Rauch
Comparison Theorem, we are going to estimate
$d^{2}\varphi(x_{0})(v)^{2}$. Consider the vector field along
$\gamma$ defined by
    $$
    Z(t)=\frac{t}{\ell}P(t),
    $$
where $P(t)$ is parallel along $\gamma$ and $P(\ell)=v$. Since $P$
is parallel (that is $P'(t)=\nabla_{\gamma'(t)}P(t)=0$), we have
    $$
    Z'(t)=\frac{1}{\ell}P(t),
    $$
and also
    $$
    \langle Z'(t), Z'(t) \rangle=\frac{1}{\ell^{2}}\|v\|_{x}^{2},
    $$
because parallel translation is a linear isometry.

On the other hand, by our assumption on the sectional curvature,
we have
\begin{eqnarray*}
& & -\langle R(\gamma'(t), P(t))\gamma'(t), P(t)\rangle\leq\\
& & K_{0}|\gamma'(t)\wedge P(t)|^{2} \leq
    K_{0} \|\gamma'(t)\|^{2} \, \|P(t)\|^{2}= K_{0}\|v\|^{2}.
\end{eqnarray*}
Therefore
\begin{eqnarray*}
& & I(Z,Z)=\int_{0}^{\ell} \langle Z'(t), Z'(t) \rangle-\langle
R(\gamma'(t), Z(t))\gamma'(t),
Z(t)\rangle dt=\\
& &\int_{0}^{\ell} \frac{1}{\ell^{2}}\|v\|_{x}^{2} -\langle
R(\gamma'(t), \frac{t}{\ell}P(t))\gamma'(t),
\frac{t}{\ell}P(t)\rangle dt=\\
& &\frac{1}{\ell^{2}}\int_{0}^{\ell}\|v\|_{x}^{2}-t^{2}\langle
R(\gamma'(t), P(t))\gamma'(t), P(t)\rangle dt\leq\\
& &
\frac{1}{\ell^{2}}\int_{0}^{\ell}\|v\|_{x}^{2}+t^{2}K_{0}\|v\|_{x}^{2}
dt=\\
& & \left(\frac{1}{\ell}+\frac{\ell}{3}K_{0} \right)\|v\|_{x}^{2}.
\end{eqnarray*}
Using Proposition \ref{Spivaks lemma} we deduce that
    $$
    I(X,X)\leq I(Z,Z)\leq\left(\frac{1}{\ell}+\frac{\ell}{3}K_{0} \right)\|v\|_{x}^{2},
    $$
hence that
    $$
    d^{2}\varphi(x)(v)^{2}=2\ell I(X,X)\leq \left(2+\frac{2\ell^{2}}{3}K_{0}
    \right)\|v\|_{x}^{2} \eqno(7).
    $$
In order to conclude the proof we only need to make sure that the
left-hand side of this inequality is nonnegative. To obtain a
lower bound for $\|d^{2}\varphi(x_{0})\|$ we will make use of the
Rauch Comparison Theorem stated above. We will compare our
manifold $M$ with a manifold $\widetilde{M}$ of constant curvature
equal to $K_{0}$, modelled in the same Hilbert space as $M$ is.

Assume first that $X$ is orthogonal to $\gamma$. Take a geodesic
$\widetilde{\gamma}:[0,\ell ]\to \widetilde{M}$ of length $\ell$,
and a vector $\widetilde{v}\in TM_{\widetilde{\gamma}(0)}$
orthogonal to $\widetilde{\gamma}'(0)$ with
$\|\widetilde{v}\|=\|X'(0)\|$. A Jacobi field $\widetilde{X}$
along $\widetilde{\gamma}$ with $\widetilde{X}(0)=0$ and
$\|\widetilde{X}'(0)\|=\|X'(0)\|$ is given by
    $$
    \widetilde{X}(t)=
    \frac{\sin\left(\sqrt{K_{0}}\, t\right)}{\sqrt{K_{0}}}
    \widetilde{P}(t),
    $$
where $\widetilde{P}(t)$ is the parallel translation of the vector
$\widetilde{v}$ along $\widetilde{\gamma}$, with
$\widetilde{P}(0)=\widetilde{v}$, see \cite[Chapter IX,
Proposition 2.12]{Lang}. Since the sectional curvature of $M$ is
less than or equal to $K_{0}$, $X(0)=0=\widetilde{X}(0)$ but
$X(t)\neq 0\neq \widetilde{X}(t)$ for $0<t\leq \ell$ (recall that
$\gamma, \widetilde{\gamma}$, having lengths less than $i(M)$, do
not have conjugate points, see \cite[Theorem 3.1 of Chapter
IX]{Lang}), $\|X'(0)\|=\|\widetilde{X}'(0)\|$, $X$ is orthogonal
to $\gamma$ and $\widetilde{X}$ is orthogonal to
$\widetilde{\gamma}$, we get from the Rauch Comparison Theorem
that
    $$
    \frac{\langle \widetilde{X}'(t), \widetilde{X}(t)\rangle}{\langle \widetilde{X}(t), \widetilde{X}(t)\rangle}
    \leq\frac{\langle X'(t), X(t)\rangle}{\langle X(t), X(t)\rangle}
    $$
for all $t\in (0, \ell]$. Since
$\|\widetilde{P}(t)\|=\|v\|=\|X(\ell)\|$ and
    $
    \widetilde{X}'(t)=\cos(\sqrt{K_{0}}\, t)\widetilde{P}(t),
    $
we deduce, by taking $t=\ell$, that
    $$
    \sqrt{K_{0}}\frac{\cos(\sqrt{K_{0}}\, \ell)}{\sin(\sqrt{K_{0}}\, \ell)}\leq
    \frac{\langle X'(\ell), X(\ell)\rangle}{\langle X(\ell), X(\ell)\rangle}=
    \frac{\langle X'(\ell), X(\ell)\rangle}{\|v\|_{x}^{2}},$$
hence
    $$
    \langle X'(\ell), X(\ell)\rangle\geq
    \sqrt{K_{0}}\frac{\cos(\sqrt{K_{0}}\, \ell)}{\sin(\sqrt{K_{0}}\,
    \ell)}\|v\|_{x}^{2}.\eqno(8)
    $$
In particular, since $\ell\leq r<\pi/2\sqrt{K_{0}}$, we obtain
that
    $$
    \langle X'(\ell), X(\ell)\rangle\geq 0.
    $$

On the other hand, if $X$ is tangent to $\gamma$ then
$$X(t)=\pm\frac{t}{\ell}\|v\|_{x}\gamma'(t),$$
and
    $$
    X'(t)=\pm\frac{1}{\ell}\|v\|_{x}\gamma'(t),$$
hence we also have
    $$
    \langle X'(\ell), X(\ell)\rangle\geq 0.
    $$

Now, from Propositions 2.3 and 2.4 of Chapter IX of \cite{Lang},
we know that every Jacobi field $X$ along $\gamma$ with $X(0)=0$
can be written in the form
    $$
    X=X^{\top}+X^{\bot},
    $$
where $X^{\top}$ and $X^{\bot}$ are Jacobi fields along $\gamma$,
$X^{\top}$ and $(X^{\top})'$ are tangent to $\gamma$, and
$X^{\bot}$ and $(X^{\bot})'$ are orthogonal to $\gamma$. In
particular $\langle X^{\top}, (X^{\bot})'\rangle=0$ and $\langle
X^{\bot}, (X^{\top})'\rangle=0$. This implies that
    $$
    \langle X'(t), X(t)\rangle=
    \langle (X^{\top})'(t), X^{\top}(t)\rangle+
    \langle (X^{\bot})'(t), X^{\bot}(t)\rangle\geq 0,
    $$
and hence
    $$
    d^{2}\varphi(x)(v)^{2}=2\ell\langle X'(\ell), X(\ell)\rangle\geq
    0, \eqno(9)
    $$
in any case.

By combining $(9)$ with equation $(7)$ above we finally get
    $$
    0\leq d^{2}\varphi(x)(v)^{2}\leq
    \left(2+\frac{2\ell^{2}}{3}K_{0}
    \right)\|v\|_{x}^{2}
    $$
for all $v\in TM_{x}$, which implies
    $$
    \|d^{2}\varphi(x)\|_{x}\leq (2+\frac{2\ell^{2}}{3}K_{0} d(x,x_{0})^{2}),
    $$
and the proof of Proposition \ref{estimation for the second
derivative of the squared distance} finishes. \hspace{6cm} $\Box$

\medskip

\begin{cor}\label{manifolds of positive injectivity radius and bounded curvature are soub}
If $i(M)>0$, $c(M)>0$, and the sectional curvature of $M$ is
bounded, then $M$ is second order uniformly bumpable.
\end{cor}
\begin{proof}
Fix any $r>0$ with $r<\min\{i(M), c(M)\}$, and define
$R=46+2K_{0}\, r^{2}$, where $K_{0}$ is a bound for the sectional
curvature of $M$. For every $\delta\in (0,r)$ find a $C^\infty$
function $\theta:\mathbb{R}\to [0,1]$ such that
\begin{enumerate}
\item[{(i)}] $\theta(t)=1$ for $t\leq 0$
\item[{(ii)}] $\theta(t)=0$ for $t\geq\delta^{2}$
\item[{(iii)}] $\|\theta'\|_{\infty}\leq 3/\delta^{2}$
\item[{(iv)}] $\|\theta''\|_{\infty}\leq 10/\delta^{4}$.
\end{enumerate}
Now, for a given $z\in M$, consider the function
$\varphi(x)=d(x,z)^{2}$, and define $b:M\to [0,1]$ by
    $$
    b (x)=
  \begin{cases}
    \theta(\varphi(x)) & \text{ if } d(x,z)\leq\delta, \\
    0 & \text{ otherwise}.
  \end{cases}
    $$
It is clear that $b$ is $C^\infty$ smooth on $M$ and $b$ satisfies
conditions $(i)$ and $(ii)$ of Definition \ref{definition of
uniform bumpability}.

In order to estimate $\|db\|_{\infty}$, first note that, since
$x\mapsto d(x,z)$ is $1$-Lipschitz, we have
    $$
    \|d\varphi(x)\|\leq 2d(x,z) \eqno(8)
    $$
for all $x\in B(z,r)$. Then, for all $x\in B(z, \delta)$,
    $$
    \|db(x)\|=|\theta'(\varphi(x))| \|d\varphi(x)\|\leq \frac{3}{\delta^{2}} 2\delta
    =\frac{6}{\delta},
    $$
and therefore
    $$
    \|db\|_{\infty}\leq\frac{6}{\delta}\leq\frac{R}{\delta},
    $$
so condition $(iii)$ of the Definition is met as well. On the
other hand, by using again $(8)$, as well as Proposition
\ref{estimation for the second derivative of the squared
distance}, we can estimate $\|d^{2}b\|_{\infty}$ as follows. For
every $x\in B(z, \delta)$ we have
\begin{eqnarray*}
& & |d^{2}b(x)(v)^{2}|=|\theta''(\varphi(x))
\left(d\varphi(x)(v)\right)^{2} +
\theta'(\varphi(x))d^{2}\varphi(x)(v)^{2}|\leq\\
& &|\theta''(\varphi(x)) \left(d\varphi(x)(v)\right)^{2}|+
|\theta'(\varphi(x))d^{2}\varphi(x)(v)^{2}|\leq\\ & &
\frac{10}{\delta^{4}}(2\delta)^{2}+\frac{3}{\delta^{2}} \left(
2+\frac{2}{3}K_{0}\, r^{2}\right)= \frac{46+2K_{0}\,
r^{2}}{\delta^{2}},
\end{eqnarray*}
which implies
    $$
    \|d^{2}\varphi\|_{\infty}\leq \frac{46+2K_{0}\,
    r^{2}}{\delta^{2}}=\frac{R}{\delta^{2}},
    $$
that is condition $(iv)$ of Definition \ref{definition of uniform
bumpability} is also satisfied.
\end{proof}

\medskip

\section{The rest of the proof}

In this section we will show that every complete Riemannian
manifold which is second-order uniformly bumpable satisfies the
natural translation to the Riemannian setting of the second-order
DGZ smooth variational principle (established in \cite{DGZ} for
all $C^2$ smooth Banach spaces). The result will be a consequence
of several auxiliary lemmas.

\begin{lem}\label{soub manifolds satisfy assumptions of superlemma}
Let $M$ be a second order uniformly bumpable Riemannian manifold.
Then there are numbers $C>1, r>0$ such that for every $p\in M$,
$\varepsilon>0$ and $\delta\in (0, r)$ there exists a $C^2$ smooth
function $b:M\to [0,\varepsilon]$ such that:
\begin{enumerate}
\item $b(p)=\varepsilon= \|b\|_{\infty}:=\sup_{x\in M}|b(x)|$.
\item $\|db\|_{\infty}:=\sup_{x\in M}\|db(x)\|_{x}\leq
C\varepsilon/\delta$.
\item $\|d^{2}b\|_{\infty}:=\sup_{x\in M}\|d^{2}b(x)\|_{x}\leq
C\varepsilon/\delta^{2}$.
\item $b(x)=0$ if $x\not\in {B(p,\delta)}$.
\end{enumerate}
In particular, $\max\{\|b\|_{\infty}, \|db\|_{\infty},
\|d^{2}b\|_{\infty}\}\leq C\varepsilon (1+1/\delta+1/\delta^{2})$.
\end{lem}
\begin{proof}
In the case when $\varepsilon=1$ we get a required $b$ from the
definition of second order uniform bumpability. If
$\varepsilon\neq 1$, it is enough to multiply $b$ by
$\varepsilon$.
\end{proof}

\medskip

\begin{lem}\label{Y is complete}
The space $Y=\{\varphi\in C^{2}(M,\mathbb{R}) \, : \, \varphi,
\|d\varphi\|, \|d^{2}\varphi\| \, \textrm{ are bounded on } \,
M\}$, endowed with the norm
    $$
    \|\varphi\|_{Y}:=\max\{\|\varphi\|_{\infty},
\|d\varphi\|_{\infty}, \|d^{2}\varphi\|_{\infty}\},
    $$
is a Banach space.
\end{lem}
\begin{proof}
The space $(Y, \|\cdot\|_{Y})$ is clearly a normed space. Let us
check that it is complete.

Let $(\varphi_{n})$ be a Cauchy sequence in $(Y, \|\cdot\|_{Y})$.
Since $\|\varphi_{n}\|_{Y}\geq \|\varphi_{n}\|_{\infty}$ and the
space of continuous bounded functions on $M$ with the norm
$\|\cdot\|_{\infty}$ is complete, we know that there exists a
continuous bounded function $\varphi:M\to\mathbb{R}$ such that
    $$
    \|\varphi_{n}-\varphi\|_{\infty}\to 0.
    $$
We have to see that $\varphi\in  C^{2}(M)$, $d\varphi$ and
$d^{2}\varphi$ are bounded, and
$\|d\varphi_{n}-d\varphi\|_{\infty}\to 0$, and
$\|d^{2}\varphi_{n}-d^{2}\varphi\|_{\infty}\to 0$. To this end,
fix a point $x\in M$ and a number $r$ with $0<r<i_{M}(x)$ and such
that both $\exp_{x}:B(0_{x}, r)\subset TM_{x}\to B(x,r)\subset M$
and its inverse $\exp_{x}^{-1}$ are $2$-Lipschitz diffeomorphisms,
and consider the functions $\psi_{n}(y)=\varphi_{n}\circ\exp_{x}$
and $\psi(y)=\varphi\circ\exp_{x}$, defined on the ball $B(0_{x},
r)$ in $TM_{x}$. We have
    $$
    d\psi_{n}(w_{y})(v)=d\varphi_{n}(y)(d\exp_{x}(w_{y})(v)),
    \eqno(1)
    $$
and $$
    d\psi(w_{y})(v)=d\varphi(y)(d\exp_{x}(w_{y})(v)), \eqno(2)
    $$
where we denote $w_{y}=\exp_{x}^{-1}(y)$. Therefore
    $$
    \sup_{w_{y}\in B(0_{x}, r)}\|d\psi_{m}(w_{y})-d\psi_{n}(w_{y})\|\leq
    2\sup_{y\in B(x, r)}\|d\varphi_{m}(y)-d\varphi_{n}(y)\|\leq
    2\|d\varphi_{m}-d\varphi_{n}\|_{\infty}
    $$
and, because $(\varphi_{n})$ is a Cauchy sequence in $Y$, the
definition of $\|\cdot\|_{Y}$ implies that the right-hand side of
the above inequality goes to $0$ as $m, n\to\infty$, which shows
that $(\psi_{n})$ is a Cauchy sequence in the space $\{ f\in
C^{1}(B(0,r)): f \textrm{ and } df \textrm{ are bounded } \}$ with
the norm $\|f\|=\max\{\|f\|_{\infty}, \|df\|_{\infty}\}$. Since
this space is complete it follows that $\psi_{n}$ converges to
some $\widetilde{\psi}\in C^{1}(B(0,r))$, and $d\psi_{n}$
converges to $d\widetilde{\psi}$, in the norm
$\|\cdot\|_{\infty}$. On the other hand we already know that
$\|\varphi_{n}-\varphi\|_{\infty}\to 0$, which implies
$\|\psi_{n}-\psi\|_{\infty}\to 0$, so $\psi=\widetilde{\psi}$ by
the uniqueness of the limit. Therefore
$\varphi=\psi\circ\exp_{x}^{-1}$ is $C^1$ on $B(x,r)$, and since
$x$ is arbitrary it follows that $\varphi\in C^{1}(M)$.

To see that $d\varphi_{n}$ converges to $d\varphi$ in the norm
$\|\cdot\|_{\infty}$, let us first observe that, equations $(1)$
and $(2)$, together with the facts that $d\exp_{x}(w_{y})$ is a
linear isomorphism and $d\psi_{n}\to d\psi$, imply that
    $$
    \|d\varphi_{n}(y)-d\varphi(y)\|_{y}\to 0, \eqno(3)
    $$
that is $d\varphi_{n}\to d\varphi$ pointwise on $TM^{*}$. Now,
since $d\varphi_{n}$ is a Cauchy sequence in the norm
$\|\cdot\|_{\infty}$, for every $\varepsilon>0$ there exists
$n_{0}\in\mathbb{N}$ such that
    $$
    \|d\varphi_{n}(y)-d\varphi_{m}(y)\|_{y}\leq\varepsilon
    \eqno(4)
    $$
for all $y\in M$, whenever $n,m\geq n_{0}$. By taking limits as
$m\to\infty$ in $(4)$, and using $(3)$ and continuity of
$\|\cdot\|_{y}$, we deduce that
    $$
    \|d\varphi_{n}(y)-d\varphi(y)\|_{y}\leq\varepsilon
    $$
for all $y\in M$, whenever $n\geq n_{0}$. This shows that
$\|d\varphi_{n}-d\varphi\|_{\infty}\to 0$.

In order to check that $\varphi\in C^{2}(M)$ and
$d^{2}\varphi_{n}\to d^{2}\varphi$, we need to use the following
Fact, which relates the second derivatives of $\psi_{n}$ and
$\varphi_{n}$.
\begin{fact}
\label{relationship through exp} Let $f:M\to\mathbb{R}$ be a $C^2$
smooth function, and define $h=f\circ\exp_{x}$ on a neighborhood
of a point $0\in TM_{x}$. Let $\widetilde{V}$ be a vector field
defined on a neighborhood of $0$ in $TM_{x}$, and consider the
vector field defined by
$V(y)=d\exp_{x}(w_{y})(\widetilde{V}(w_{y}))$ on a neighborhood of
$x$, where $w_{y}:=\exp_{x}^{-1}(y)$, and let
    $$
    \sigma_{y}(t)=\exp_{x}(w_{y}+t\widetilde{V}(w_{y})).
    $$
Then we have that
    $$
    D^{2}h(\widetilde{V},\widetilde{V})(w_{y})=D^{2}f(V,V)(y)+\langle\nabla f(y), \sigma_{y}''(0)\rangle.
    $$
\end{fact}
\noindent The proof of this fact is just a calculation, see
\cite[Lemma 2.7]{AFS2}. We apply this with $\widetilde{V}(w)=v$, a
constant field on $TM_{x}$, $f=\varphi_{n}$, $h=\psi_{n}$ to
obtain
    $$
d^{2}\psi_{n}(w_{y})(v,v)=d^{2}\varphi_{n}(y)(d\exp_{x}(w_{y})(v),
d\exp_{x}(w_{y})(v))+\langle\nabla \varphi_{n}(y),
\sigma_{y}''(0)\rangle, \eqno(5)
    $$
where
    $$
    \sigma_{y}(t)=\exp_{x}(w_{y}+tv).
    $$
Since the vector field $M\ni y\mapsto \sigma_{y}''(0)\in TM$ is
continuous and $\sigma_{x}''(0)=0$, we can assume without loss of
generality that $r>0$ is small enough so that
$\|\sigma_{y}''(0)\|_{y}\leq 1$ for all $y\in B(x,r)$. Then, again
using the fact that $\exp_{x}$ is $2$-biLipschitz on $B(x,r)$, we
get from $(5)$ that
\begin{eqnarray*}
& &\sup_{w_{y}\in B(0, r)}\|d^{2}\psi_{n}(w_{y})-d^{2}\psi_{m}(w_{y})\|\leq\\
& & 2\sup_{y\in B(x,r)}\|d^{2}\varphi_{n}(y)-d^{2}\varphi_{n}(y)\|
    + \sup_{y\in B(x,r)}\|d\varphi_{m}(y)-d\varphi_{n}(y)\|\leq\\
& & 2\|d^{2}\varphi_{n}-d^{2}\varphi_{m}\|_{\infty} +
\|d\varphi_{n}-d\varphi_{m}\|_{\infty},
\end{eqnarray*}
which implies that $(\psi_{n})$ is a Cauchy sequence in the space
$\{ f\in C^{2}(B(0,r)): f, \, df,  \textrm{ and } d^{2}f \textrm{
are bounded } \}$ with the norm $\|f\|=\max\{\|f\|_{\infty},
\|df\|_{\infty}, \|d^{2}f\|_{\infty}\}$. This space is well known
to be complete, hence there exists some $\widetilde{\psi}\in
C^{2}(B(0,r))$ such that $\psi_{n}$ converges to
$\widetilde{\psi}$, and $d^{2}\psi_{n}$ converges to
$d^{2}\widetilde{\psi}$, in the norm $\|\cdot\|_{\infty}$. Since
we already know that $\|\psi_{n}-\psi\|_{\infty}\to 0$, we get
that $\psi=\widetilde{\psi}$, and therefore
$\varphi=\psi\circ\exp_{x}^{-1}$ is $C^2$ on $B(x,r)$. It follows
that $\varphi\in C^{2}(M)$.

Moreover, equation $(5)$ (and the same equation replacing
$\varphi_{n}$ with $\varphi$ and $\psi_{n}$ with $\psi$), together
with the facts that $d\exp_{x}(w_{y})$ is a linear isomorphism,
and that $d\varphi_{n}\to\varphi$, imply that
    $$
    \|d^{2}\varphi_{n}(y)-d^{2}\varphi(y)\|_{y}\to 0,
    $$
for each $y\in B(x,r)$, that is $d^{2}\varphi_{n}\to d^{2}\varphi$
pointwise on $T_{2,s}(M)$.

By combining this with the fact that $(d^{2}\varphi_{n})$ is a
Cauchy sequence in the norm $\|\cdot\|_{\infty}$, one can easily
deduce (as in the case of $(d\varphi_{n})$) that
$\|d^{2}\varphi_{n}-d^{2}\varphi\|_{\infty}\to 0$.
\end{proof}

\medskip

In the sequel $B(\varphi, r)$ stands for the open ball of center
$\varphi$ and radius $r$ in the Banach space $Y$.
\begin{lem}\label{the superlemma}
Let $M$ be a complete metric space, and $(Y,\|\cdot\|)$ be a
Banach space of real-valued bounded and continuous functions on
$M$ satisfying the following conditions:
\begin{enumerate}
\item $\|\varphi\|\geq\|\varphi\|_{\infty}=\sup\{|\varphi(x)|:x\in
M\}$ for every $\varphi\in Y$.
\item There are numbers $C>1, r>0$ such that for every $p\in M$,
$\varepsilon>0$ and $\delta\in (0, r)$ there exists a function
$b\in Y$ such that $b(p)=\varepsilon$, $\|b\|_{Y}\leq C\varepsilon
(1+1/\delta+1/\delta^{2})$, and $b(x)=0$ if $x\not\in
{B(p,\delta)}$.
\end{enumerate}
Let $f:M\rightarrow \mathbb{R}\cup\{+\infty\}$ be a lower
semicontinuous function which is bounded below and such that
$Dom(f)=\{x\in M | f(x)<+\infty\}\neq\emptyset$. Then, the set $G$
of all the functions $\varphi\in Y$ such that $f+\varphi$ attains
a strong minimum in $M$ contains a $G_{\delta}$ dense subset of
$Y$.
\end{lem}
\begin{proof} The proof follows the lines of that of Lemma 3.13 in
\cite{AFL2} (which in turn is very similar to the original proof
of \cite{DGZsvp}), with some small changes. We write the complete
proof (rather than just indicating the changes) for the reader's
convenience and for completeness.

Take a number $N\in\mathbb{N}$ such that $N\geq 1/r$, and for
every $n\in\mathbb{N}$ with $n\geq N$, consider the set
$$U_n=\{\varphi\in Y | \, \exists \, x_0 \in M \,
:\, (f+\varphi)(x_0)<\inf\{(f+\varphi)(x)|x\in M\backslash
B(x_{0},\frac{1}{n})\}\}.$$
\begin{claim}
$U_n$ is open.
\end{claim}
\noindent Take $\varphi\in U_n$. By the definition of $U_n$ there
exists $x_{0}\in M$ such that
$(f+\varphi)(x_0)<\inf\{(f+\varphi)(x)|x\in M\backslash
B(x_{0},\frac{1}{n})\}$. Set $2\rho=\inf\{(f+\varphi)(x)|x\in
M\backslash B(x_{0},\frac{1}{n})\}-(f+\varphi)(x_0)>0$. Then,
since $\|\cdot\|_Y\geq \|\cdot\|_{\infty}$, we get that
$B_{Y}(\varphi, \rho)\subset B_{\infty}(\varphi,\rho)\subset U_n$.
\begin{claim}
$U_n$ is dense in $Y$.
\end{claim}
\noindent Take $\varphi\in Y$ and $\varepsilon>0$. Since
$f+\varphi$ is bounded below there exists $x_0 \in M$ such that
$(f+\varphi)(x_0)< \inf\{(f+\varphi)(x)|x\in M\}+\varepsilon$. Set
now $\delta=1/n<r$, and use condition (2) to find a function $b\in
Y$ such that $b(x_0)=\varepsilon$, $\|b\|_{Y}\leq
C(n^{2}+n+1)\varepsilon$, and $b(x)=0$ for $x\not\in
B(x_{0},\frac{1}{n})$. Then
$(f+\varphi)(x_0)-b(x_0)<\inf\{(f+\varphi)(x)|x\in M\}$ and, if we
define $h=-b$, we have
    $$(f+\varphi+h)(x_0)<\inf\{(f+\varphi)(x)|x\in M\}\leq
    \inf\{(f+\varphi)(x)|x\not\in
B(x_0,\frac{1}{n})\}.$$ Since $\inf\{(f+\varphi)(x)|x\not\in
B(x_0,\frac{1}{n})\}=\inf\{(f+\varphi+h)(x)|x\not\in
B(x_0,\frac{1}{n})\}$, it is obvious that the above inequality
implies that $\varphi+h\in U_n$. On the other hand, we have
$\|h\|_{Y}\leq C(n^{2}+n+1)\varepsilon$. Since $C$ and $n$ are
fixed and $\varepsilon$ can be taken to be arbitrarily small, this
shows that $\varphi\in\overline{U_{n}}$, and $U_{n}$ is dense in
$Y$.

\medskip

Therefore we can apply Baire's theorem to conclude that the set
$G=\bigcap_{n=N}^{\infty} U_n$ is a $G_\delta$ dense subset of
$Y$. Now we must show that if $\varphi\in G$ then $f+\varphi$
attains a strong minimum in $M$. For each $n\geq N$, take $x_n \in
M$ such that $(f+\varphi)(x_n)<\inf\{(f+\varphi)(x)|x\not\in
B(x_n,\frac{1}{n})\}$. Clearly, $x_k\in B(x_n,\frac{1}{n})$ if
$k\geq n$, which implies that $(x_n)_{n=N}^{\infty}$ is a Cauchy
sequence in $M$ and therefore converges to some $x_0 \in M$. Since
$f$ is lower semicontinuous and $\bigcap_{n=N}^{\infty}B(x_{0},
1/n)=\{x_{0}\}$, we get
\begin{eqnarray*}
&&(f+\varphi)(x_{0})\leq\liminf(f+\varphi)(x_n)\leq
\liminf[\inf\{(f+\varphi)(x)|x\in M\setminus
B(x_0,\frac{1}{n})\}]\\ &&=\inf\{\inf\{(f+\varphi)(x)|x\in
M\setminus B(x_{0},\frac{1}{n})\}
: n\in\mathbb{N}, n\geq N\}\\
&&=\inf\{(f+\varphi)(x)|x\in M\setminus \{x_0\}\},
\end{eqnarray*}
which means that $f+\varphi$ attains a global minimum at $x_0\in
M$.

Finally, let us check that in fact $f+\varphi$ attains a strong
minimum at the point $x_{0}$. Suppose $\{y_n\}$ is a sequence in
$M$ such that $(f+g)(y_n)\rightarrow (f+g)(x_0)$ and $(y_n)$ does
not converge to $x_0$. We may assume $d(y_n, x_0)\geq \varepsilon$
for all $n$. Bearing in mind this inequality and the fact that
$x_{0}=\lim x_{n}$, we can take $k\in\mathbb{N}$ such that $d(x_k,
y_n)>\frac{1}{k}$ for all $n$, and therefore
    $$(f+\varphi)(x_0)\leq
(f+\varphi)(x_k)<\inf\{(f+\varphi)(x) | x\notin B(x_{k},
\frac{1}{k}) \}\leq(f+\varphi)(y_n)$$ for all $n$, which
contradicts the fact that $(f+\varphi)(y_n)\rightarrow
(f+\varphi)(x_0)$.
\end{proof}

\medskip

By combining Lemmas \ref{soub manifolds satisfy assumptions of
superlemma}, \ref{Y is complete} and \ref{the superlemma} we
immediately deduce the following.
\begin{thm}\label{soub manifolds satisfy the sosvp}
Let $M$ be a second order uniformly bumpable complete Riemannian
manifold. Then, for every lower semicontinuous function $f:M\to
(-\infty, \infty]$ which is bounded below, with
$f\not\equiv+\infty$, and for every $\varepsilon>0$, there exists
a $C^2$ smooth function $\varphi: M\to\mathbb{R}$ such that
\begin{enumerate}
\item $f-\varphi$ attains its strong minimum on $M$
\item $\|\varphi\|_{\infty}<\varepsilon$
\item $\|d\varphi\|_{\infty}<\varepsilon$
\item $\|d^{2}\varphi\|_{\infty}<\varepsilon$
\end{enumerate}
\end{thm}

\medskip

Finally, Theorem \ref{main theorem} follows from Theorem \ref{soub
manifolds satisfy the sosvp} and Corollary \ref{manifolds of
positive injectivity radius and bounded curvature are soub}.



\end{document}